\documentclass[12pt]{amsart}
\usepackage{amsmath, amssymb, amsthm, url, scrextend, mathtools}
\usepackage[all]{xy}
\usepackage{color}
\usepackage[top=2.5cm, bottom=2.5cm, left=2.5cm, right=2.5cm]{geometry}

\newcommand{\nc}{\newcommand}
\nc{\lei}{\le^\oo}
\nc{\zero}{\mathbf{0}}
\nc{\zar}{\mathfrak{zar}}
\nc{\card}[1]{\left|#1\right|}
\nc{\bbN}{\mathbb{N}}
\nc{\beq}{\begin{eqnarray*}}\nc{\eeq}{\end{eqnarray*}}
\nc{\mbq}{\mb{?}}
\nc{\mb}[1]{{\mbox{\textbf{#1}}}}
\nc{\nop}{$\times$}
\nc{\fbn}{\!\!\fbox{\!\nop\!}\!\!}
\nc{\yup}{\checkmark}
\nc{\forces}{\Vdash}
\nc{\name}[1]{\dot{#1}}
\nc{\eseq}[1]{#1_1, \allowbreak #1_2, \allowbreak\dotsc}
\nc{\tf}{\my{FINISHED THUS FAR}}
\nc{\FU}{Fr\'echet--Urysohn}
\nc{\gs}{$\gamma$~space}
\nc{\Ga}{\Gamma}
\nc{\Om}{\Omega}
\nc{\cozOm}{\Omega_{\mathrm{coz}}}
\nc{\ctblOm}{\Omega_{\mathrm{ctbl}}}
\nc{\smallbinom}[2]{{\begin{psmallmatrix} #1\\ #2 \end{psmallmatrix}}}
\nc{\two}{\{0,1\}}
\nc{\productive}[2]{\bigl(#1,\allowbreak #2\bigr)^\x}
\nc{\Sel}{\mathsf{S}}
\nc{\sset}[2]{{\{\,#1 : #2\,\}}}
\nc{\smb}[1]{{\!\!\mb{#1}\!\!}}
\nc{\cZ}{\mathcal{Z}}
\nc{\medset}[2]{{\biggl\{\,#1 : #2\,\biggr\}}}
\nc{\smallmedset}[2]{{\bigl\{\,#1 : #2\,\bigr\}}}
\nc{\set}[2]{{\left\{\,#1 : #2\,\right\}}}
\nc{\seq}[2]{{\la\, #1 : #2\,\ra}}
\nc{\cube}{(\Cantor)^\bbN}
\nc{\Match}{\op{Match}}
\nc{\concat}[1]{\hat{\phantom{a}}\langle #1\rangle}
\nc{\poset}{\mathbb{P}}
\nc{\fn}[1]{{\op{Fn}(#1\times\w,2)}}
\nc{\linadd}{\op{linadd}}
\nc{\nonprod}{\non^\x}
\nc{\alephes}{{\aleph_0}}
\nc{\my}[1]{\textcolor{red}{#1}}
\nc{\Cp}{\op{C}}
\nc{\Cpp}{\op{C_p}}
\nc{\Bp}{\op{B}_p}
\nc{\Pa}[8]{\bibitem{#1} {#2}, \emph{#3}, {#4} \textbf{#5} ({#6}), {#7}--{#8}.}
\nc{\tPa}[5]{\bibitem{#1} {#2}, \emph{#3}, {#4}, to appear.}
\nc{\sPa}[4]{\bibitem{#1} {#2}, \emph{#3}, {#4}, submitted.}
\nc{\Bc}[9]{\bibitem{#1} {#2}, \emph{#3}, in: \textbf{#4} (#5), #6 #7, #8--#9.}
\nc{\fD}{\mathfrak{D}}
\nc{\fX}{\mathfrak{X}}
\nc{\Onbd}{\Op_{\mathrm{nbd}}} 
\nc{\Omnb}{\Om_{\mathrm{nbd}}} 
\nc{\od}{\mathfrak{od}}
\nc{\Setting}[7]{\xymatrix@R=4pt@C=7pt{#1\ar@{-}[r]&#2\ar@{-}[r]&#3\\&#4\ar@{-}[u]\\
#5\ar@{-}[uu]\ar@{-}[r] & #6\ar@{-}[u]\ar@{-}[r] & #7\ar@{-}[uu]}}
\nc{\mx}[1]{\begin{matrix}#1\end{matrix}}
\nc{\plim}{p\txt{-}\lim}
\nc{\Bgp}{{\Z^\bbN}}
\nc{\Cgp}{{{\Z_2}^\bbN}}
\nc{\Cite}[1]{\textbf{[#1]}}
\nc{\Next}[1]{{#1^+}}
\nc{\Fr}{\mathit{F\!r}}
\nc{\intvl}[2]{{[#1(#2),\allowbreak #1(#2\!+\!1))}}
\nc{\Bdd}{\mathbf{B}}
\nc{\Dfin}{\mathfrak{D}_\mathrm{fin}}
\nc{\grbl}{{\mbox{\textit{\tiny gp}}}}
\nc{\bbP}{\mathbb{P}}
\nc{\BOfat}{\B_{\Om_{\mathrm{fat}}}}
\nc{\Bgood}{\B_{\mathrm{good}}}
\nc{\compactN}{\cl{\mathbb{N}}}
\nc{\blocks}[2]{\op{cl}_{#2}(#1)}
\nc{\blocksplus}[2]{\op{cl}^+_{#2}(#1)}
\nc{\arx}[1]{\texttt{http://arxiv.org/math/#1}}
\nc{\bq}{\begin{quote}}
\nc{\eq}{\end{quote}}
\nc{\cl}[1]{\overline{#1}}
\nc{\CH}{the Continuum Hypothesis}
\nc{\MA}{Martin's Axiom}
\nc{\Bfat}{\B_\mathrm{fat}}
\nc{\inv}{^{-1}}
\nc{\Cantor}{{\two^\bbN}}
\nc{\bP}{\mathbf{P}}
\nc{\bA}{\mathbf{A}}
\nc{\bB}{\mathbf{B}}
\nc{\bof}{\op{\fb}}
\nc{\bofF}{\bof(\cF)}
\nc{\sr}[3]{\underset{\mbox{#3}}{\mbox{#1}}}
\nc{\gig}{\gimel}
\nc{\gns}{\sone(\Om,\gig)}
\nc{\nsr}[2]{#1}
\nc{\Srg}{{\mathbb{S}}}
\nc{\Srgs}{{\mathbb{S}^*}}
\nc{\NN}{{\bbN^{\bbN}}}
\nc{\ZN}{{\Z^{\bbN}}}
\nc{\NNup}{{\bbN^{\uparrow\bbN}}}
\nc{\Pof}{\op{P}}
\nc{\PN}{{\Pof(\bbN)}}
\nc{\roth}{{[\bbN]^{\mbox{\tiny $\infty$}}}} 
\nc{\Fin}{\mathrm{Fin}}
\nc{\ici}{[\bbN]^{ \infty, \infty}}
\nc{\Inc}{{\compactN^{\uparrow\bbN}}}
\nc{\powInc}[1]{{\big(\Inc\big)^{#1}}}
\nc{\powFin}[1]{{\big(\Fin\big)^{#1}}}
\nc{\powPN}[1]{{\big(\PN\big)^{#1}}}
\nc{\NcompactN}{{\compactN^\bbN}}
\nc{\Uarrow}{\smash{\big\uparrow}}
\nc{\LE}{\preccurlyeq}
\nc{\GE}{\succcurlyeq}
\nc{\op}{\operatorname}
\nc{\im}{\op{im}}
\nc{\Span}{\op{span}}
\nc{\maxfin}{\op{maxfin}}
\nc{\ran}{\op{range}}
\nc{\iso}{\cong}
\nc{\Madd}{{\M}^\star}
\nc{\cI}{\mathcal{I}}
\nc{\cJ}{\mathcal{J}}
\nc{\scrA}{\mathscr{A}}
\nc{\scrB}{\mathscr{B}}
\nc{\scrC}{\mathscr{C}}
\nc{\scrD}{\mathscr{D}}
\nc{\scrF}{\mathscr{F}}
\nc{\scrK}{\mathscr{K}}
\nc{\A}{\forall}
\nc{\B}{\mathrm{B}}
\nc{\cB}{\mathcal{B}}
\nc{\BS}{\mathbf{B}(\mathcal{S})}
\nc{\BF}{\mathbf{B}(\mathcal{F})}
\nc{\BU}{\mathbf{B}(\mathcal{U})}
\nc{\cSp}{\mathcal{S}^+}
\nc{\cFp}{\mathcal{F}^+}
\nc{\cUp}{\mathcal{U}^+}

\nc{\BG}{\Ga_\mathrm{Bor}}
\nc{\BO}{\Om_\mathrm{Bor}}

\nc{\BL}{\B_\Lambda}
\nc{\BT}{\B_\Tau}
\nc{\BTstar}{\B_{\Tau^*}}
\nc{\DO}{\cD_\Om}
\nc{\KO}{\cK_\Om}
\nc{\CG}{C_\Ga}
\nc{\CL}{C_\Lambda}
\nc{\CT}{C_\Tau}
\nc{\CTstar}{C_{\Tau^*}}
\nc{\CO}{C_\Om}
\nc{\sfC}{\mathsf{C}}
\nc{\sfD}{\mathsf{D}}
\nc{\bD}{\mathbf{D}}
\nc{\Tau}{\mathrm{T}}
\nc{\cA}{\mathcal{A}}
\nc{\cK}{\mathcal{K}}
\nc{\cD}{\mathcal{D}}
\nc{\cF}{\mathcal{F}}
\nc{\cS}{\mathcal{S}}
\nc{\cT}{\mathcal{T}}
\nc{\cG}{\mathcal{G}}
\nc{\cY}{\mathcal{Y}}
\nc{\J}{\mathcal{J}}
\nc{\cL}{\mathcal{L}}
\nc{\cM}{\mathcal{M}}
\nc{\cN}{\mathcal{N}}
\nc{\cE}{\mathcal{E}}
\nc{\cH}{\mathcal{H}}
\nc{\cO}{\mathcal{O}}
\nc{\Op}{\mathrm{O}}
\nc{\rmA}{\mathrm{A}}
\nc{\rmB}{\mathrm{B}}
\nc{\rmD}{\mathrm{D}}
\nc{\rmP}{\mathrm{P}}
\nc{\cC}{\mathcal{C}}
\nc{\cP}{\mathcal{P}}
\nc{\bbQ}{\mathbb{Q}}
\nc{\bbR}{\mathbb{R}}
\nc{\bbZ}{\mathbb{Z}}
\nc{\cU}{\mathcal{U}}
\nc{\Un}{\bigcup}
\nc{\cV}{\mathcal{V}}
\nc{\cW}{\mathcal{W}}
\nc{\Z}{{\mathbb Z}}
\nc{\Impl}{\Rightarrow}
\long\def\forget#1\forgotten{}
\nc{\ft}{\mathfrak{t}}
\nc{\fb}{\mathfrak{b}}
\nc{\fc}{\mathfrak{c}}
\nc{\fd}{\mathfrak{d}}
\nc{\fg}{\mathfrak{g}}
\nc{\oo}{\infty}
\nc{\fr}{\mathfrak{r}}
\nc{\fk}{\mathfrak{k}}
\nc{\bidi}{\mathfrak{bidi}}
\nc{\fu}{\mathfrak{u}}
\nc{\fh}{\mathfrak{h}}
\nc{\fq}{\mathfrak{q}}
\nc{\fp}{\mathfrak{p}}
\nc{\fj}{\mathfrak{j}}
\nc{\fs}{\mathfrak{s}}
\nc{\w}{\omega}
\nc{\x}{\times}
\nc{\Iff}{\Leftrightarrow}

\nc{\nin}{\notin}
\nc{\cat}{\hat{\ }}
\nc{\sub}{\subseteq}
\nc{\spst}{\supseteq}
\nc{\sm}{\setminus}
\nc{\as}{\subseteq^*}
\nc{\aspst}{\prescript{*}{}{\spst}\ }
\nc{\les}{\le^*}
\nc{\leinf}{\le^{\infty}}
\nc{\leS}{\le_{\mathcal{S}}}
\nc{\leF}{\le_{\mathcal{F}}}
\nc{\leU}{\le_{\mathcal{U}}}
\nc{\rest}{\restriction}

\nc{\la}{\langle}
\nc{\ra}{\rangle}
\nc{\E}{\exists}
\nc{\dom}{\op{dom}}
\nc{\cov}{\op{cov}}
\nc{\add}{\op{add}}
\nc{\cof}{\op{cof}}
\nc{\cf}{\op{cf}}
\nc{\non}{\op{non}}
\nc{\unif}{\op{non}}
\nc{\COV}{\op{COV}}
\nc{\ADD}{\op{ADD}}
\nc{\COF}{\op{COF}}
\nc{\NON}{\op{NON}}
\nc{\bfP}{\mathbf{P}}
\nc{\impl}{\to}
\nc{\Lp}{\mathcal{L_\p}}
\nc{\Wlog}{without loss of generality}
\newtheorem{thm}{Theorem}
\nc{\bthm}{\begin{thm}} \nc{\ethm}{\end{thm}}
\newtheorem{prop}[thm]{Proposition}
\nc{\bprp}{\begin{prop}} \nc{\eprp}{\end{prop}}
\newtheorem{fact}[thm]{Fact}
\nc{\bfct}{\begin{fact}} \nc{\efct}{\end{fact}}
\newtheorem{prob}[thm]{Problem}
\nc{\bprb}{\begin{prob}} \nc{\eprb}{\end{prob}}
\newtheorem{lem}[thm]{Lemma}
\nc{\blem}{\begin{lem}} \nc{\elem}{\end{lem}}
\newtheorem{claim}[thm]{Claim}
\nc{\bclm}{\begin{claim}} \nc{\eclm}{\end{claim}}
\newtheorem{cor}[thm]{Corollary}
\nc{\bcor}{\begin{cor}} \nc{\ecor}{\end{cor}}
\newtheorem{conj}[thm]{Conjecture}
\nc{\bcnj}{\begin{conj}} \nc{\ecnj}{\end{conj}}
\theoremstyle{definition}
\newtheorem{defn}[thm]{Definition}
\nc{\bdfn}{\begin{defn}} \nc{\edfn}{\end{defn}}
\newtheorem{obs}[thm]{Observation}
\nc{\bobs}{\begin{obs}} \nc{\eobs}{\end{obs}}
\theoremstyle{remark}
\newtheorem{rem}[thm]{Remark}
\nc{\brem}{\begin{rem}} \nc{\erem}{\end{rem}}
\newtheorem{cnv}[thm]{Convention}
\nc{\bcnv}{\begin{cnv}} \nc{\ecnv}{\end{cnv}}
\newtheorem{exam}[thm]{Example}
\nc{\bexm}{\begin{exam}} \nc{\eexm}{\end{exam}}
\nc{\bpf}{\begin{proof}} \nc{\epf}{\end{proof}}
\nc{\be}{\begin{enumerate}}
\nc{\ee}{\end{enumerate}}
\nc{\bi}{\begin{itemize}}
\nc{\bimy}{\my{\begin{itemize}}
\nc{\eimy}{\end{itemize}}}
\nc{\itm}{\item}
\nc{\ei}{\end{itemize}}
\nc{\ed}{\end{document}}

\nc{\czOm}{\Om_{\mathrm{cz}}}

\nc{\cozga}{\smallbinom{\cozOm}{\Ga}}
\nc{\GNga}{{\smallbinom{\Om}{\Ga}}}
\nc{\ctblga}{\smallbinom{\ctblOm}{\Ga}}
\nc{\borga}{\smallbinom{\BO}{\Ga}}
\nc{\czga}{\smallbinom{\czOm}{\Ga}}

\DeclareMathOperator{\Int}{Int}
\nc{\sone}{\mathsf{S}_1}
\nc{\sfin}{\mathsf{S}_\mathrm{fin}}
\nc{\ufin}{\mathsf{U}_\mathrm{fin}}
\nc{\Split}{\mathsf{Split}}

\nc{\gone}{\mathsf{G}_1}    \nc{\gfin}{\mathsf{G}_\mathrm{fin}}

\title[Strongly sequentially separable spaces]{Strongly sequentially separable function spaces,\\
	 via selection principles}

\author[A. Osipov]{Alexander V.~Osipov}
\address{Alexander V.~Osipov, Krasovskii Institute of Mathematics and Mechanics, Ural Federal University, Ural State University of Economics, 620219, Yekaterinburg, Russia}
\email{OAB@list.ru}

\author[P. Szewczak]{Piotr Szewczak}
\address{Piotr Szewczak, Institute of Mathematics, Faculty of Mathematics and Natural Science College of Sciences, Cardinal Stefan Wyszy\'nski University in Warsaw, Warsaw, Poland}
\email{p.szewczak@wp.pl}
\urladdr{www.piotrszewczak.pl}

\author[B. Tsaban]{Boaz Tsaban}
\address{Boaz Tsaban, Department of Mathematics, Bar-Ilan University, Ramat Gan, Israel}
\email{tsaban@math.biu.ac.il}
\urladdr{http://math.biu.ac.il/~tsaban}

\subjclass[2010]{
37F20, 
26A03, 
03E75, 
54C35. 
}
		
\keywords{strong sequential separability, function spaces, selection principles, Gerlits--Nagy, $\gamma$-property, Borel function, $\GNga$, $\borga$, $\gamma$-set, $\Cpp$-space.}

\begin{document}

\begin{abstract}
A separable space is \emph{strongly sequentially separable} if, for each
countable dense set, every point in the space is a limit of a sequence
from the dense set.
We consider this and related properties, for the spaces of
continous and Borel real-valued functions on Tychonoff spaces, with the topology
of pointwise convergence.
Our results solve a problem stated by Gartside, Lo, and Marsh.
\end{abstract}

\maketitle

\section{Introduction}

We apply methods of selection principles to a problem of Gartside, Lo, and Marsh~\cite[Problem~19]{glm}.

By \emph{space} we mean a Tychonoff topological space.
A space is \emph{Fr\'echet--Urysohn} if each point in the closure of a set
is a limit of a sequence from the set.
A separable space is \emph{strongly sequentially separable (SSS)}~\cite{Mat00} if, for each
countable dense set, every point in the space is a limit of a sequence from
the dense set.
Every separable Fr\'echet--Urysohn space is strongly sequentially separable,
but not conversely~\cite[Example~2.4]{bbm13}.

For a space $X$, let $\Cp(X)$ and $\B(X)$ be the spaces of continuous and Borel, respectively,
real-valued functions on $X$, with the topology of pointwise convergence.
We are only concerned with uncountable spaces. In this case, the space $\B(X)$ is never Fr\'echet--Urysohn. Indeed, for an uncountable space, the constant function $\mathbf{1}$ is in the closure of the set of characteristic functions of finite subsets of the space and there is no sequence in the set converging to $\mathbf{1}$.
Strong sequential separability is hereditary for separable dense subspaces.
Thus, if the space $\Cp(X)$ is separable, we have the following implications.
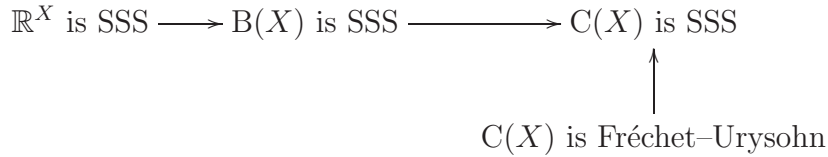
\[
\xymatrix{
\bbR^X\text{ is SSS}\ar[r] & \B(X)\text{ is SSS}\ar[r] & \Cp(X)\text{ is SSS}\\
& & \Cp(X)\text{ is Fr\'echet--Urysohn}\ar[u]
}
\]

It is consistent that the properties in this diagram hold only
for countable spaces $X$ and are, thus, equivalent~\cite[Corollary~17]{glm}.
This motivates the following problem~\cite[Problem~19]{glm}.

\bprb[{Gartside--Lo--Marsh}]\label{prb:glm}
Is there, consistently, a space $X$ such that the space $\Cp(X)$ is strongly sequentially separable but not Fr\'echet--Urysohn, and the space $\bbR^X$ is not strongly sequentially separable?
\eprb

We solve this problem, and all other problems suggested by the above diagram.
To this end, we extend Arhangel'ski\u{\i}'s local-to-global duality, dualize these problems
to ones concerning covering properties, and apply the theory of selection principles.

\section{Local-to-global Duality}

A \emph{cover} of a space is a family of \emph{proper} subsets whose union is the entire space.
For families $\bA$ and $\bB$ of covers of a space, the property that every
cover in the family $\bA$ has a subcover in the family $\bB$ is denoted $\smallbinom{\bA}{\bB}$.
An \emph{$\w$-cover} is a cover such that each finite subset of the space is contained in some set from the cover.
A \emph{$\gamma$-cover} is an infinite cover such that each point of the space belongs to all but finitely many sets from the cover.

An \emph{open cover} is a cover by open sets. Similarly, we define \emph{Borel cover}, \emph{clopen cover},
etc.
Given a space, let $\Om$, $\ctblOm$, $\cozOm$, $\BO$ and $\Ga$, be the families of \emph{open} $\omega$-covers, \emph{countable open} $\omega$-covers, \emph{countable cozero} $\omega$-covers, \emph{countable Borel} $\omega$-covers, and $\gamma$-covers, respectively.

The property $\GNga$ is the celebrated $\gamma$-property of Gerlits and Nagy, who proved that a space has this property if and only if the space $\Cp(X)$ is Fr\'echet--Urysohn~\cite[Theorem~2]{gn}.




\blem\label{lem:Bsssequiv} Let $X$ be a space with a coarser second countable topology.
The following assertions are equivalent:
\be
\item The space $\B(X)$ is strongly sequentially separable.
\item The space $X$ has the property $\borga$.
\ee
\elem
\bpf
(1) $\Impl$ (2): Since the space $X$ has a coarser second countable topology, there is a countable dense set $H$ in the space $\Cp(X)$ ~\cite[Theorem~1]{nobl} and hence $H$ is dense in the space $\B(X)$. Let $\cU\in\BO(X)$.
For a Borel set $U\sub X$ and a function $h\in H$,
let $f_{U,h}\in\B(X)$ be the function such that $f_{U,h}\restriction U:=h\restriction U$
and $f_{U,h}\restriction (X\sm U):=1$.
The set $D:=\sset{ f_{U,h}}{U\in\cU, h\in H }$ is a countable dense subset of $\B(X)$.
By (1), there is a sequence
$\sset{f_{U_n,h_n}}{n\in \bbN}$ in the set $D$, converging to the zero function $\bf{0}$.
Let $F$ be a finite subset of $X$.
The set $W:=\sset{f\in \B(X)}{f[F]\sub(-1,1)}$ is a neighborhood of $\bf{0}$ in $\B(X)$.
For a natural number $n$, if $f_{U_n,h_n}\in W$, then $F\sub U_n$.
Since all but finitely many elements of the sequence belong to the set $W$, we have $\sset{U_n}{n\in\bbN}\in\Ga(X)$.
Thus, the space $X$ satisfies $\borga$.

(2) $\Impl$ (1): The property $\borga$ implies that every point in the closure of a countable
set in $\B(X)$ is the limit of a sequence from that set~\cite[Lemma~2.8]{orts}.
\epf

Let $\bbN$ be the set of natural numbers.
For infinite sets $a,b\sub\bbN$ we write $a\as b$ if the set $a\sm b$ is finite.
A \emph{pseudointersection} of a family of infinite
sets is an infinite set $a$ with $a\as b$ for all sets $b$ in the family.
A subfamily of $\roth$ is \emph{centered} if the finite intersections of its elements, are infinite.
Let $\fp$ be the minimal cardinality of a family of infinite subsets of $\bbN$
that is centered and has no pseudointersection.
The hypothesis $\aleph_1<\fp$ and its negation ($\aleph_1=\fp$) are both consistent~\cite[Theorem~5.1]{vD}.
Information about the cardinal number $\fp$ is available, for example, in van Douwen's
survey~\cite{vD}.

Gartside, Lo and Marsh proved that a Tychonoff product $\bbR^X$
is strongly sequentially separable if and only if $\card{X} < \fp$~\cite[Theorem~11]{glm}.
They also proved that
a function space $\Cp(X)$ is strongly sequentially separable if and only if
the space $X$ has a coarser second countable topology,
and every coarser second countable topology for $X$ satisfies $\GNga$~\cite[Theorem~16]{glm}.
The property $\GNga$ implies $\smallbinom{\cozOm}{\Ga}$. Bonanzinga, Cammaroto and Matveev proved that a space $X$ has the property $\smallbinom{\cozOm}{\Ga}$ if and only if every coarser second countable topology for the space $X$ has the property $\smallbinom{\Om}{\Ga}$~\cite[Theorem~54]{bcm}.
In summary, for spaces $X$ with a coarser second countable topology,
the diagram from the previous section dualizes to the following one.
\[
\xymatrix{
\card{X}<\fp\ar[r] & X\text{ satisfies }\binom{\BO}{\Ga}\ar[r] & X\text{ satisfies }\binom{\cozOm}{\Ga}\\
& & X\text{ satisfies }\binom{\Om}{\Ga}\ar[u]
}
\]
Problem~\ref{prb:glm} is thus reduced to the following problem.

\bprb
\label{prb:dual}
Is there, consistently, a space $X$ with a coarser second countable topology,
that satisfies $\smallbinom{\cozOm}{\Ga}$ but not $\smallbinom{\Om}{\Ga}$, with $\card{X}\geq\fp$?
\eprb

We will solve this problem, as well as its variations.

\section{The problems and their solutions}

Consider the positions in the diagrams from the previous section.
Write there ``$\bullet$'' if the property holds, and ``$\circ$'' if it does not.
For example, sets $X\sub\bbR$ of cardinality smaller than $\fp$ realize the following setting.
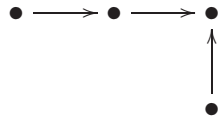
\[
\xymatrix{
\bullet \ar[r] &	\bullet \ar[r] & \bullet\\
&	& \bullet\ar[u]
}
\]
that will be denoted
$\begin{smallmatrix}
\bullet & \bullet & \bullet\\& & \bullet
\end{smallmatrix}$
for brevity.
We consider the consistency of all settings that are not ruled out by the implications
in the diagram. These are the following settings:

\medskip
\[
\begin{smallmatrix}
\circ & \circ & \circ\\& & \circ
\end{smallmatrix}
\qquad
\begin{smallmatrix}
\circ & \circ & \bullet\\& & \circ
\end{smallmatrix}
\qquad
\begin{smallmatrix}
\circ & \bullet & \bullet\\& & \circ
\end{smallmatrix}
\qquad
\begin{smallmatrix}
\circ & \circ & \bullet\\& & \bullet
\end{smallmatrix}
\qquad
\begin{smallmatrix}
\bullet & \bullet & \bullet\\& & \circ
\end{smallmatrix}
\qquad
\begin{smallmatrix}
\circ & \bullet & \bullet\\& & \bullet
\end{smallmatrix}
\qquad
\begin{smallmatrix}
\bullet & \bullet & \bullet\\& & \bullet
\end{smallmatrix}
\]
\medskip

Problem~\ref{prb:glm} asks whether either of the the settings
$\begin{smallmatrix}
	\circ & \bullet & \bullet\\& & \circ
\end{smallmatrix}$
or
$\begin{smallmatrix}
	\circ & \circ & \bullet\\& & \bullet
\end{smallmatrix}$
is consistent.

The following proposition is a variation of an earlier result~\cite[Example~4.7]{ospyt}.

\bprp[$\begin{smallmatrix}
	\bullet & \bullet & \bullet\\& & \circ
\end{smallmatrix}$]
The following assertions are equivalent:
\be
\item There is a space $X$ such that the space $\bbR^X$ is strongly sequentially separable,
but the space $\Cp(X)$ is not Fr\'echet--Urysohn.
\item $\aleph_1<\fp$.
\ee
\eprp
\bpf
Recall that the space $\bbR^X$ is strongly sequentially separable if and only if $\card{X}<\fp$.

$(1)\Impl(2)$: The given space $X$ has $\card{X}<\fp$. Had it been countable, the space
$\Cp(X)$ would have been metrizable.

$(2)\Impl(1)$:
A discrete space of cardinality $\aleph_1$ is not Lindel\"of, and thus
does not satisfy $\GNga$. Apply duality.
\epf

The following folklore fact implies that discrete spaces of cardinality $\fp$ or greater have
none of the studied properties.

\blem[$\begin{smallmatrix}
	\circ & \circ & \circ\\& & \circ
\end{smallmatrix}$]
\label{lem:nonGNctbl} Let $X$ be a set. Then $\card{X}<\fp$ if and only if every countable $\omega$-cover consisting of subsets of $X$ contains a $\gamma$-cover.
\elem
\bpf
($\Rightarrow$) Let $\sset{U_n}{n\in\bbN}$ be a countable $\omega$-cover consisting of subsets of $X$.
For each element $x\in X$, let $a_x:=\sset{n\in\bbN}{x\in U_n}$, an infinite subset of $\bbN$.
Since  $\card{X}<\fp$, the family $\sset{a_x}{x\in X}$ has a pseudointersection $a$.
Then $\sset{U_n}{n\in a}$ is a $\gamma$-cover of $X$.

($\Leftarrow$) Assume that $\card{X}\ge\fp$. We may assume that $X\sub \roth$ where $\roth$ is the family of infinite subsets of $\bbN$. Then $X$ is a family of infinite subsets of $\bbN$ of cardinality $\fp$,
that is centered and has no pseudointersection.
The family of sets $\sset{U_n}{n\in\bbN}$, defined by $U_n:=\sset{x\in X}{n\in x}$ for natural numbers $n$,
is a countable $\omega$-cover of $X$ and has no subfamily in $\Ga$.
\epf

A theorem of Galvin and Miller~\cite[Theorem~2]{gami} asserts that, if $\fp=\card{\bbR}$, then
there is a set $X\sub\bbR$ of cardinality $\fp$, satisfying $\GNga$.
The Galvin--Miller Theorem is refined by Theorem~\ref{thm:tower} of Orenshtein and Tsaban~\cite[Theorem~3.6]{ot}.
Since this result is central
to the remainder of this paper, we include here a simpler proof, due to the third named author~\cite{tslec}.

We identify the Cantor space $\Cantor$ with the family $\PN$ of all subsets of the set $\bbN$.
Thus, we view the space $\PN$ as a subset of the real line.
The space $\PN$ splits into two subspaces: the family of infinite subsets of $\bbN$, denoted $\roth$, and the family of finite subsets of $\bbN$, denoted $\Fin$.
We identify every set $a\in\roth$ with its increasing enumeration, an element of the Baire space $\NN$.
Thus, for a natural number $n$, $a(n)$ is the $n$-th element in the increasing enumeration of the set $a$.
This way, we have $\roth\sub\NN$, and the topology of the space $\roth$ (a subspace
of the Cantor space $\PN$) coincides
with the subspace topology induced by $\NN$. When an element of $\roth$ is viewed as an element of $\NN$,
we refer to it as a \emph{function}.

For functions $a,b\in\roth$, we write $a\les b$ if the set $\sset{n}{b(n)<a(n)}$ is finite.
Let $A\sub\roth$.
For a function $b\in \roth$, we write $A\les b$ if $a\les b$ for all functions $a\in A$.
The set $A$ is \emph{unbounded} if there is no function $b\in \roth$ with $A\les b$.
Let $\fb$ be the minimal cardinality of an unbounded set in $\roth$.
A set $\sset{x_\alpha}{\alpha<\fb}\sub\roth$ is an \emph{unbounded tower} if it is unbounded and for all ordinal numbers $\alpha,\beta<\fb$ with $\alpha<\beta$, we have $x_\alpha \aspst   x_\beta$.
An unbounded tower of cardinality $\fp$ exists if (and only if) $\fp=\fb$~\cite[Lemma~3.3]{ot}.

\bthm[Orenshtein--Tsaban]\label{thm:tower}
For each unbounded tower $T\sub\roth$ of cardinality $\fp$,
the set $T\cup\Fin$ of real numbers satisfies $\GNga$.
\ethm

In order to prove Theorem~\ref{thm:tower}, we need the following notions and auxiliary results.
Let $n,m$ be natural numbers with $n<m$.
Define $(n,m):=\sset{i\in\bbN}{n<i<m}$.
A set $a\in \roth$ \emph{omits} the interval $(n,m)$ if $a\cap (n,m)=\emptyset$.
For a space $X$, let $\Om(X)$ be  the family of all \emph{open $\omega$-covers} of $X$,
 and $\Ga(X)$ be the family of all \emph{open $\gamma$-covers} of $X$.

\blem[{Galvin--Miller\cite[Lemma~1.2]{gami}}]\label{lem:GM}
Let $\cU$ be a family of open sets in $\PN$ such that $\cU\in\Om(\Fin)$.
There are a function $b\in\roth$ and distinct 
sets $\eseq{U}\in \cU$ such that for each element
$x\in \roth$ and all natural numbers~$n$:
\[
\text{If }x\cap((b(n), b(n+1))=\emptyset,\text{ then }x\in U_n.
\]
\elem

\blem[{Folklore~\cite[Lemma~2.13]{MHP}}]\label{lem:unbdd}
Let $Y$ be a subset of $\roth$.
The set $Y$ is unbounded if and only if, for each function $b\in\roth$, there is a set $a\in Y$ that omits infinitely many intervals $(b(n),b(n+1))$.
\elem

\blem\label{lem:key}
Let $X\sub\PN$ be a set such that $\Fin\sub X$ and $\card{X}<\fp$.
Let $\cU$ be a family of open sets in $\PN$ such that $\cU\in\Om(X)$, and  $Y$ be an unbounded set in $\roth$.
There are a set $a\in Y$, and sets $\eseq{U}\in\cU$ such that $\sset{U_n}{n\in\bbN}\in \Ga(X)$, and for each element
 $x\in\roth$ and all natural numbers $n$:
\[\text{If }x\sm\{\,1,\dotsc,n\,\}\sub a,\text{ then }x\in\bigcap_{k\geq n}U_k.
\]
\elem

\bpf
Since $\card{X}<\fp$, the set $X$ satisfies $\GNga$~\cite[Proposition~2]{reclaw}.
Let $\cV\in\Ga(X)$ be a subfamily of $\cU$.
By Lemma~\ref{lem:GM}, there are a function $b\in\roth$, and distinct sets $\eseq{V}\in\cV$ such that for each element $x\in \roth$, and all natural numbers $i$:
\begin{equation}\label{eq:step1}
	\text{If }x\cap(b(i),b(i+1))=\emptyset,\text{ then }x\in V_i.
\end{equation}
By Lemma~\ref{lem:unbdd}, there is a set $a\in Y$ such that the set
\[
c:=\set{i\in\bbN}{a\cap(b(i),b(i+1))=\emptyset}
\]
is infinite.
Fix a natural number $n$.
Let $k$ be a natural number with $n\leq k$, and $x\in\roth$ be an element such that $x\sm\{\,1,\dotsc,n\,\}\sub a$.
Then $n\leq c(k)$, and we have
\[
x\cap (b(c(k)),b(c(k)+1))\sub a\cap (b(c(k)),b(c(k)+1))=\emptyset.
\]
By~\eqref{eq:step1}, we have $x\in V_{c(k)}$.
Thus, $x\in \bigcap_{k\geq n}V_{c(k)}$.

Since $\cV\in\Ga(X)$, we have $\sset{V_{c(i)}}{i\in\bbN}\in\Ga(X)$.
\epf

\begin{proof}[Proof of Theorem~\ref{thm:tower}]
Let $\sset{x_\alpha}{\alpha<\fb}\sub\roth$ be an unbounded tower.
Let $X:=\Fin\cup\sset{x_\alpha}{\alpha<\fb}$, and for ordinal numbers $\gamma<\fb$, let
$X_\gamma:=\Fin\cup\sset{x_\alpha}{\alpha<\gamma}$.
Let $\cU\in\Om(X)$. Fix an ordinal number $\gamma_0<\fb$.
By induction, for a natural number $m>0$, we proceed as follows.
By Lemma~\ref{lem:key}, there are an ordinal number $\gamma_m<\fb$, and a subfamily  $\sset{U^{(m)}_n}{n\in\bbN}\in\Ga(X_{\gamma_{m-1}})$ of $\cU$ such that, for each element $x\in\roth$ and all natural numbers $n$:
	
\begin{equation}\label{eq:key}
\text{If }x\sm\{\,1,\dotsc,n\,\}\sub x_{\gamma_m},\text{ then }x\in \bigcap_{k\geq n}U^{(m)}_k.
\end{equation}

Let $\gamma:=\sup_n\gamma_n$.
There is a function $g\in \roth$ such that
$x_\gamma\sm\{\,1,\dotsc,g(n)\,\}\sub x_{\gamma_n}$ for all natural numbers $n$.
Fix an ordinal number $\alpha$ with $\gamma\leq\alpha<\fb$.
Since $x_\alpha\as x_\gamma$, we have
\[
x_\alpha\sm\{\,1,\dotsc,g(n)\,\}\sub x_\gamma\sm\{\,1,\dotsc,g(n)\,\}\sub x_{\gamma_n},
\]
for all but finitely many natural numbers $n$.
By~\eqref{eq:key}, we have $x_\alpha\in \bigcap_{k\geq g(n)}U^{(n)}_{k}$ for all but finitely many natural numbers $n$.
Thus, for any function $h\in\roth$ with $g\les h$, we have $\sset{U^{(n)}_{h(n)}}{n\in\bbN}\in\Ga(\sset{x_\alpha}{\gamma\leq\alpha<\fb})$.
	
For each element $x\in X_\gamma$, and each natural number $n$, define
\[
f_x(n):=\min\medset{m\in\bbN}{x\in\bigcap_{k\geq m}U^{(n)}_k}
\]
if the set is nonempty, and $f_x(n):=0$ otherwise.
Since $\card{X_\gamma}<\fb$, there is a function $h\in\roth$ such that $\sset{f_x}{x\in X_\gamma}\cup \{g\}\les h$, and the sets $U^{(n)}_{h(n)}$ are distinct.
Then $\sset{U^{(n)}_{h(n)}}{n\in\bbN}\in\Ga(X_\gamma)$.
Since $\sset{U^{(n)}_{h(n)}}{n\in\bbN}\in\Ga(\sset{x_\alpha}{\gamma\leq\alpha<\fb})$ as well,
we have  $\sset{U^{(n)}_{h(n)}}{n\in\bbN}\in\Ga(X)$.
\epf
	
\section{Subsets of the Real, Michael, and Sorgenfrey line}\label{sec:exm}

The \emph{Michael line}~\cite{michael} is the set $\PN$, with the topology where the points of the set $\roth$ are isolated, and the neighborhoods of the points of the set $\Fin$ are those induced by the Cantor space topology on $\PN$.
The \emph{Sorgenfrey line}~\cite{srg} is the set $\bbR$ with the topology generated by the
half-open intervals $[a,b)$, for $a,b\in\bbR$.

The forthcoming Theorem~\ref{thm:exm}(2) solves the problem of Gartside--Lo--Marsh Problem (Problem~\ref{prb:glm}).
Recall that an unbounded tower in $\roth$ of cardinality $\fp$ exists if and only if $\fp=\fb$.
It is consistent that $\aleph_1<\fp=\fb$, e.g., it holds assuming the Martin Axiom with the negation of \CH{}.

\bthm\label{thm:exm} Let $T\sub\roth$ be an  unbounded tower  of cardinality $\fp$.
\be
\item
$\begin{psmallmatrix}
	\circ & \circ & \bullet\\& & \bullet
\end{psmallmatrix}$
As a subset of $\bbR$, the set $T\cup\Fin$ satisfies $\GNga$ but not $\borga$.
\item
$\begin{psmallmatrix}
	\circ & \circ & \bullet\\& & \circ
\end{psmallmatrix}$
Assume that $\aleph_1<\fp$. As a subset of the Michael line,
the set $T\cup\Fin$ satisfies $\ctblga$ but neither $\GNga$ nor  $\borga$.
\ee
\ethm
\bpf
(1)
By Theorem~\ref{thm:tower}, the set $T\cup\Fin$ satisfies $\GNga$.
The set $T$ is centered and has no pseudointersection.
Thus, the set $T$ does not satisfy $\borga$ ~\cite[Lemma~24]{scts}.
Since the set $T$ is a Borel subset of $T\cup\Fin$, and the property
$\borga$ is hereditary for Borel subsets ~\cite[Theorem~48]{scts}, the set $T\cup\Fin$ does not satisfy $\borga$, too.

(2)
For a set $U\sub\PN$, let $\Int(U)$ be the interior of the set $U$ in the Cantor space topology on $\PN$.
If $\cU\in\Om(\Fin)$ is a family of open sets in the Michael line, then $\sset{\Int(U)}{U\in\cU}\in \Om(\Fin)$.
Thus, Lemma~\ref{lem:GM} holds for a family of open sets in the Michael line.
By Lemma~\ref{lem:nonGNctbl}, every space of cardinality smaller than $\fp$ satisfies $\ctblga$.
Thus, Lemma~\ref{lem:key} holds for a countable family of open sets in the Michael line.
Consequently, the proof of Theorem~\ref{thm:tower} actually establishes that the set $T\cup\Fin$,
as a subspace of the Michael line, satisfies $\ctblga$.

Write $T=\sset{x_\alpha}{\alpha<\fb}$ with $x_\alpha\as x_\beta$ for $\beta<\alpha$.
The set $A:=\sset{x\in T}{x_{\w_1}\as x}$ has cardinality $\aleph_1$.
The set $A$ is $\mathrm{F}_\sigma$ in the Cantor space topology and, in particular, in the
Michael line topology.
Thus, the space $T\cup\Fin$ has an uncountable discrete $F_\sigma$ subset.
Since the Lindel\"of property is hereditary for $\mathrm{F}_\sigma$ subsets, the space $T\cup\Fin$ is not Lindel\"of.
Every space with the property $\GNga$ is Lindel\"of.
Thus, the space $T\cup\Fin$ does not satisfy $\GNga$.

By (1), since every Borel set in the Cantor space is also Borel in the Michael line,
the space $T\cup\Fin$ does not satisfy $\borga$.
\epf

\bcor
Let $T\sub\roth$ be an  unbounded tower  of cardinality $\fp$.
\be
\item
For the real line topology, the space $\Cp(T\cup\Fin)$ is Fr\'echet-Urysohn but
the space $\B(T\cup\Fin)$ is not strongly sequentially separable.
\item
Assume that $\aleph_1<\fp$. For the Michael line topology,
the space $\Cp(T\cup\Fin)$ is strongly sequentially separable and not Fr\'echet--Urysohn,
and the space $\B(T\cup\Fin)$ is not strongly sequentially separable.\qed
\ee
\ecor

Assuming \CH{}, there is an uncountable set of real numbers satisfying $\borga$~(\cite[Theorem~4.1]{br},~\cite[Theorem~5]{mil1}).
If $\aleph_1<\fp$, then any subset of real numbers of cardinality $\aleph_1$ satisfies $\borga$~\cite[Lemma~22, Theorem~27(1)]{scts}.

\bthm
Let $X\sub\bbR$ be an uncountable set satisfying $\borga$.
\be
\item
As a subset of $\bbR$, the set $X$ satisfies $\GNga$.
\item
As a subset of the Sorgenfrey line, the set $-X\cup X$ satisfies $\borga$ but not $\GNga$.
\ee
In particular, if \CH{} holds, we obtain the setting
$\begin{smallmatrix}
\circ & \bullet & \bullet\\& & \bullet
\end{smallmatrix}$
from (1), and the setting
$\begin{smallmatrix}
\circ & \bullet & \bullet\\& & \circ
\end{smallmatrix}$
from (2).
If $\aleph_1<\fp$,
we obtain the settings
$\begin{smallmatrix}
\bullet & \bullet & \bullet\\& & \bullet
\end{smallmatrix}$
and
$\begin{smallmatrix}
\bullet & \bullet & \bullet\\& & \circ
\end{smallmatrix}$.
\ethm
\bpf
(1) For subsets of $\bbR$, the property $\borga$ implies $\GNga$.

(2) Let $Y\sub\bbR$ be an uncountable set satisfying $\borga$.
The disjoint union $Y\sqcup Y$ satisfies $\borga$ as well: Let $\cU$ be a countable Borel $\omega$-cover of
$Y\sqcup Y$. The family
\[
\cV:=\sset{U\cap V}{U\sqcup V\in\cU,\ U\sub Y\sqcup\emptyset, V\sub \emptyset\sqcup Y}
\]
is a  countable Borel $\omega$-cover of $Y$. Let $\cW\sub\cV$ be a $\gamma$-cover of $Y$.
Then the family
\[
\sset{U\sqcup V\in\cU}{U\cap V\in\cW, \ U\sub Y\sqcup\emptyset, V\sub \emptyset\sqcup Y}
\]
is a $\gamma$-cover of $Y\sqcup Y$.

The set $X:=Y\cup \sset{-y}{y\in Y}$, a continuous image of the space $Y\sqcup Y$, satisfies $\borga$, too.
Consider this set as a subspace of the Sorgenfrey line.
Since the Borel sets in the real line and the Sorgenfrey line are the same, the space $X$ satisfies $\borga$.

The product space $X\x X$ contains the uncountable closed discrete set $\sset{(x,-x)}{x\in X}$, and thus
does not satisfy $\GNga$.
The property $\GNga$ is preserved by finite powers~\cite[Theorem~3.6]{coc2}.
Thus, the space $X$ does not satisfy $\GNga$.
\epf

\bcor
Let $X\sub\bbR$ be an uncountable set satisfying $\borga$.
As a subset of the Sorgenfrey line, the space $\Cp(-X\cup X)$ is not Fr\'echet--Urysohn,
but the space $\B(-X\cup X)$ is strongly sequentially separable.\qed
\ecor

\section{Additional results}

The space from Theorem~\ref{thm:exm}(2) has the  property $\ctblga$, that is formally stronger than $\cozga$.
In the forthcoming Proposition~\ref{prp:GNctblfoGNctbl}, we show that the properties $\ctblga$ and $\cozga$ are different.

A family of sets is \emph{almost disjoint} if the intersection of any pair of sets of this family is finite.
For an almost disjoint family $A$ in $\roth$, the \emph{Mr\'owka--Isbell} space $\Psi(A)$ is the set $A\cup\bbN$, with the points of $\bbN$ isolated, and with
the sets $\{a\}\cup a\sm b$ (for $b\in\Fin$) as  neighborhoods of the points $a\in A$.

\bprp\label{prp:GNctblfoGNctbl}
There is a maximal almost disjoint family $A$ in $\roth$ such that the Mr\'owka--Isbell space $\Psi(A)$ satisfies $\cozga$ but not $\ctblga$.
\eprp

\bpf
There is a maximal almost disjoint family $A$ in $\roth$, of cardinality $\card{\bbR}$, such that the space $\Psi(A)$ satisfies $\cozga$~\cite[Example~61 and Theorem~54]{bcm}.
Let $A=\sset{a_r}{r\in\bbR}$.
Since $\bbR$ does not satisfy $\ctblga$, there is a family $\cU\in\ctblOm(\bbR)$ with no subfamily in $\Ga(\bbR)$.
For each set $U\in\cU$, let $U':=\sset{a_r}{r\in U}\cup\bbN$.
The family $\sset{U'}{U\in\cU}$ is in $\ctblOm(\Psi(A))$ and has no subfamily in $\Ga(\Psi(A))$.
Thus, the space $\Psi(A)$ does not satisfy $\ctblga$.
\epf

A space is \emph{projectively $\GNga$} if each continuous second countable image of this space satisfies $\GNga$~\cite{bcm}.

\bprp\label{prp:sssequiv}
For a space $X$, the following assertions are
equivalent:
\be
\item The space $\Cp(X)$ is strongly sequentially separable.
\item The space $X$ has a coarser second countable topology, and it is projectively $\GNga$.
\ee
\eprp

\bpf
$(1)\Impl(2)$:
By a result of Noble~\cite[Theorem~1]{nobl}, the space $X$ has a coarser second countable topology.
In order to prove that the space $X$ is projectively $\GNga$, we show that it satisfies the equivalent property $\cozga$~\cite[Theorem~54]{bcm}.
Let $F\sub \Cp(X)$ be a countable set such that the family $\cU=\sset{ f\inv[\bbR\sm\{0\}]}{f\in F}$ is an $\w$-cover of $X$.
Let $\cB$ be a countable basis of $\bbR$, and $\cB'$ be a countable basis of a coarser topology on $X$.
Let $Y$ be the set $X$ with the topology generated by the family $\sset{f\inv[B]}{B\in \cB}\cup\cB'$.
The space $Y$ is second countable.
By a result of Gartside, Lo, and Marsh~\cite[Theorem~16]{glm}, the space $Y$ satisfies $\GNga$.
Since $\cU\in \Om(Y)$, the family $\cU$ contains a cover $\cV\in\Ga(Y)$.
Thus, $\cV\in\Ga(X)$.

$(2)\Impl(1)$:
By~(2), every coarser second countable topology for the space $X$ satisfies $\GNga$.
By a result of Gartside, Lo and Marsh~\cite[Theorem~16]{glm}, the space $\Cp(X)$ is strongly sequentially separable.
\epf

\subsection*{Acknowledgments} We thank the referee for useful comments and corrections.

\ed